\documentstyle[12pt]{article}

\newcommand{\tlowername}[2]%
{$\stackrel{\makebox[1pt]{#1}}%
{\begin{picture}(0,0)%
\put(0,0){\makebox(0,6)[t]{\makebox[1pt]{$#2$}}}%
\end{picture}}$}%

\newcommand{\AR}[1]%
{\begin{picture}(#1,0)%
\put(0,0){\vector(1,0){#1}}%
\end{picture}}%

\newcommand{\DOTAR}[1]%
{\NUMBEROFDOTS=#1%
\divide\NUMBEROFDOTS by 3%
\begin{picture}(#1,0)%
\multiput(0,0)(3,0){\NUMBEROFDOTS}{\circle*{1}}%
\put(#1,0){\vector(1,0){0}}%
\end{picture}}%

\newcommand{\MONO}[1]%
{\begin{picture}(#1,0)%
\put(0,0){\vector(1,0){#1}}%
\put(2,-2){\line(0,1){4}}%
\end{picture}}%

\newcommand{\EPI}[1]%
{\begin{picture}(#1,0)(-#1,0)%
\put(-#1,0){\vector(1,0){#1}}%
\put(-6,-2){\line(0,1){4}}%
\end{picture}}%

\newcommand{\BIMO}[1]%
{\begin{picture}(#1,0)(-#1,0)%
\put(-#1,0){\vector(1,0){#1}}%
\put(-6,-2){\line(0,1){4}}%
\put(-#1,-2){\hspace{2pt}\line(0,1){4}}%
\end{picture}}%

\newcommand{\BIAR}[1]%
{\begin{picture}(#1,4)%
\put(0,0){\vector(1,0){#1}}%
\put(0,4){\vector(1,0){#1}}%
\end{picture}}%

\newcommand{\EQL}[1]%
{\begin{picture}(#1,0)%
\put(0,1){\line(1,0){#1}}%
\put(0,-1){\line(1,0){#1}}%
\end{picture}}%

\newcommand{\ADJAR}[1]%
{\begin{picture}(#1,4)%
\put(0,0){\vector(1,0){#1}}%
\put(#1,4){\vector(-1,0){#1}}%
\end{picture}}%

\newcommand{\BKAR}[1]%
{\begin{picture}(#1,0)%
\put(#1,0){\vector(-1,0){#1}}%
\end{picture}}%

\newcommand{\BKDOTAR}[1]%
{\NUMBEROFDOTS=#1%
\divide\NUMBEROFDOTS by 3%
\begin{picture}(#1,0)%
\multiput(#1,0)(-3,0){\NUMBEROFDOTS}{\circle*{1}}%
\put(0,0){\vector(-1,0){0}}%
\end{picture}}%

\newcommand{\BKMONO}[1]%
{\begin{picture}(#1,0)(-#1,0)%
\put(0,0){\vector(-1,0){#1}}%
\put(-2,-2){\line(0,1){4}}%
\end{picture}}%

\newcommand{\BKEPI}[1]%
{\begin{picture}(#1,0)%
\put(#1,0){\vector(-1,0){#1}}%
\put(6,-2){\line(0,1){4}}%
\end{picture}}%

\newcommand{\BKBIMO}[1]%
{\begin{picture}(#1,0)%
\put(#1,0){\vector(-1,0){#1}}%
\put(6,-2){\line(0,1){4}}%
\put(#1,-2){\hspace{-2pt}\line(0,1){4}}%
\end{picture}}%

\newcommand{\BKBIAR}[1]%
{\begin{picture}(#1,4)%
\put(#1,0){\vector(-1,0){#1}}%
\put(#1,4){\vector(-1,0){#1}}%
\end{picture}}%

\newcommand{\BKADJAR}[1]%
{\begin{picture}(#1,4)%
\put(0,4){\vector(1,0){#1}}%
\put(#1,0){\vector(-1,0){#1}}%
\end{picture}}%

\newcommand{\lowername}[2]%
{$\stackrel{\makebox[1pt]{#1}}%
{\begin{picture}(0,0)%
\truex{600}%
\put(0,0){\makebox(0,\value{x})[t]{\makebox[1pt]{$#2$}}}%
\end{picture}}$}%

\newcommand{\hcase}[2]%
{\makebox[0pt]%
{\raisebox{-1pt}[0pt][0pt]{#1{#2}}}}%

\newcommand{\Hcase}[3]%
{\makebox[0pt]
{\raisebox{-1pt}[0pt][0pt]%
{$\stackrel{\makebox[0pt]{$\textstyle{#2}$}}{#1{#3}}$}}}%

\newcommand{\hcasE}[3]%
{\makebox[0pt]%
{\raisebox{-9pt}[0pt][0pt]%
{\lowername{#1{#3}}{#2}}}}%

\newcommand{\hbicase}[2]%
{\makebox[0pt]%
{\raisebox{-2.5pt}[0pt][0pt]{#1{#2}}}}%

\newcommand{\Hbicase}[4]%
{\makebox[0pt]
{\raisebox{-10.5pt}[0pt][0pt]%
{$\stackrel{\makebox[0pt]{$\textstyle{#2}$}}%
{\mbox{\lowername{#1{#4}}{#3}}}$}}}%

\newcommand{\EAR}[1]%
{\begin{picture}(#1,0)%
\put(0,0){\vector(1,0){#1}}%
\end{picture}}%

\newcommand{\EDOTAR}[1]%
{\truex{100}\truey{300}%
\NUMBEROFDOTS=#1%
\divide\NUMBEROFDOTS by \value{y}%
\begin{picture}(#1,0)%
\multiput(0,0)(\value{y},0){\NUMBEROFDOTS}%
{\circle*{\value{x}}}%
\put(#1,0){\vector(1,0){0}}%
\end{picture}}%

\newcommand{\EMONO}[1]%
{\begin{picture}(#1,0)%
\put(0,0){\vector(1,0){#1}}%
\truex{300}\truey{600}%
\put(\value{x},-\value{x}){\line(0,1){\value{y}}}%
\end{picture}}%

\newcommand{\EEPI}[1]%
{\begin{picture}(#1,0)(-#1,0)%
\put(-#1,0){\vector(1,0){#1}}%
\truex{300}\truey{600}\truez{800}%
\put(-\value{z},-\value{x}){\line(0,1){\value{y}}}%
\end{picture}}%

\newcommand{\EBIMO}[1]%
{\begin{picture}(#1,0)(-#1,0)%
\put(-#1,0){\vector(1,0){#1}}%
\truex{300}\truey{600}\truez{800}%
\put(-\value{z},-\value{x}){\line(0,1){\value{y}}}%
\put(-#1,-\value{x}){\hspace{3pt}\line(0,1){\value{y}}}%
\end{picture}}%

\newcommand{\EBIAR}[1]%
{\truex{400}%
\begin{picture}(#1,\value{x})%
\put(0,0){\vector(1,0){#1}}%
\put(0,\value{x}){\vector(1,0){#1}}%
\end{picture}}%

\newcommand{\EEQL}[1]%
{\begin{picture}(#1,0)%
\truex{200}%
\put(0,\value{x}){\line(1,0){#1}}%
\put(0,0){\line(1,0){#1}}%
\end{picture}}%

\newcommand{\EADJAR}[1]%
{\truex{400}%
\begin{picture}(#1,\value{x})%
\put(0,0){\vector(1,0){#1}}%
\put(#1,\value{x}){\vector(-1,0){#1}}%
\end{picture}}%

\newcommand{\ear}%
{\hspace{\SOURCE\unitlength}%
\hcase{\EAR}{\ARROWLENGTH}}%

\newcommand{\Ear}[1]%
{\hspace{\SOURCE\unitlength}%
\Hcase{\EAR}{#1}{\ARROWLENGTH}}%

\newcommand{\eaR}[1]%
{\hspace{\SOURCE\unitlength}%
\hcasE{\EAR}{#1}{\ARROWLENGTH}}%

\newcommand{\edotar}%
{\hspace{\SOURCE\unitlength}%
\hcase{\EDOTAR}{\ARROWLENGTH}}%

\newcommand{\Edotar}[1]%
{\hspace{\SOURCE\unitlength}%
\Hcase{\EDOTAR}{#1}{\ARROWLENGTH}}%

\newcommand{\edotaR}[1]%
{\hspace{\SOURCE\unitlength}%
\hcasE{\EDOTAR}{#1}{\ARROWLENGTH}}%

\newcommand{\emono}%
{\hspace{\SOURCE\unitlength}%
\hcase{\EMONO}{\ARROWLENGTH}}%

\newcommand{\Emono}[1]%
{\hspace{\SOURCE\unitlength}%
\Hcase{\EMONO}{#1}{\ARROWLENGTH}}%

\newcommand{\emonO}[1]%
{\hspace{\SOURCE\unitlength}%
\hcasE{\EMONO}{#1}{\ARROWLENGTH}}%

\newcommand{\eepi}%
{\hspace{\SOURCE\unitlength}%
\hcase{\EEPI}{\ARROWLENGTH}}%

\newcommand{\Eepi}[1]%
{\hspace{\SOURCE\unitlength}%
\Hcase{\EEPI}{#1}{\ARROWLENGTH}}%

\newcommand{\eepI}[1]%
{\hspace{\SOURCE\unitlength}%
\hcasE{\EEPI}{#1}{\ARROWLENGTH}}%

\newcommand{\ebimo}%
{\hspace{\SOURCE\unitlength}%
\hcase{\EBIMO}{\ARROWLENGTH}}%

\newcommand{\Ebimo}[1]%
{\hspace{\SOURCE\unitlength}%
\Hcase{\EBIMO}{#1}{\ARROWLENGTH}}%

\newcommand{\ebimO}[1]%
{\hspace{\SOURCE\unitlength}%
\hcasE{\EBIMO}{#1}{\ARROWLENGTH}}%

\newcommand{\eiso}%
{\hspace{\SOURCE\unitlength}%
\Hcase{\EAR}{\cong}{\ARROWLENGTH}}%

\newcommand{\Eiso}[1]%
{\hspace{\SOURCE\unitlength}%
\Hcase{\EAR}{\cong#1}{\ARROWLENGTH}}%

\newcommand{\eisO}[1]%
{\hspace{\SOURCE\unitlength}%
\hcasE{\EAR}{\cong#1}{\ARROWLENGTH}}%

\newcommand{\ebiar}%
{\hspace{\SOURCE\unitlength}%
\hbicase{\EBIAR}{\ARROWLENGTH}}%

\newcommand{\Ebiar}[2]%
{\hspace{\SOURCE\unitlength}%
\Hbicase{\EBIAR}{#1}{#2}{\ARROWLENGTH}}%

\newcommand{\eeql}%
{\hspace{\SOURCE\unitlength}%
\hbicase{\EEQL}{\ARROWLENGTH}}%

\newcommand{\eadjar}%
{\hspace{\SOURCE\unitlength}%
\hbicase{\EADJAR}{\ARROWLENGTH}}%

\newcommand{\Eadjar}[2]%
{\hspace{\SOURCE\unitlength}%
\Hbicase{\EADJAR}{#1}{#2}{\ARROWLENGTH}}%

\newcommand{\WAR}[1]%
{\begin{picture}(#1,0)%
\put(#1,0){\vector(-1,0){#1}}%
\end{picture}}%

\newcommand{\WDOTAR}[1]%
{\truex{100}\truey{300}%
\NUMBEROFDOTS=#1%
\divide\NUMBEROFDOTS by \value{y}%
\begin{picture}(#1,0)%
\multiput(#1,0)(-\value{y},0){\NUMBEROFDOTS}%
{\circle*{\value{x}}}%
\put(0,0){\vector(-1,0){0}}%
\end{picture}}%

\newcommand{\WMONO}[1]%
{\begin{picture}(#1,0)(-#1,0)%
\put(0,0){\vector(-1,0){#1}}%
\truex{300}\truey{600}%
\put(-\value{x},-\value{x}){\line(0,1){\value{y}}}%
\end{picture}}%

\newcommand{\WEPI}[1]%
{\begin{picture}(#1,0)%
\put(#1,0){\vector(-1,0){#1}}%
\truex{300}\truey{600}\truez{800}%
\put(\value{z},-\value{x}){\line(0,1){\value{y}}}%
\end{picture}}%

\newcommand{\WBIMO}[1]%
{\begin{picture}(#1,0)%
\put(#1,0){\vector(-1,0){#1}}%
\truex{300}\truey{600}\truez{800}%
\put(\value{z},-\value{x}){\line(0,1){\value{y}}}%
\put(#1,-\value{x}){\hspace{-3pt}\line(0,1){\value{y}}}%
\end{picture}}%

\newcommand{\WBIAR}[1]%
{\truex{400}%
\begin{picture}(#1,\value{x})%
\put(#1,0){\vector(-1,0){#1}}%
\put(#1,\value{x}){\vector(-1,0){#1}}%
\end{picture}}%

\newcommand{\WADJAR}[1]%
{\truex{400}%
\begin{picture}(#1,\value{x})%
\put(0,\value{x}){\vector(1,0){#1}}%
\put(#1,0){\vector(-1,0){#1}}%
\end{picture}}%

\newcommand{\war}%
{\hspace{\SOURCE\unitlength}%
\hcase{\WAR}{\ARROWLENGTH}}%

\newcommand{\War}[1]%
{\hspace{\SOURCE\unitlength}%
\Hcase{\WAR}{#1}{\ARROWLENGTH}}%

\newcommand{\waR}[1]%
{\hspace{\SOURCE\unitlength}%
\hcasE{\WAR}{#1}{\ARROWLENGTH}}%

\newcommand{\wdotar}%
{\hspace{\SOURCE\unitlength}%
\hcase{\WDOTAR}{\ARROWLENGTH}}%

\newcommand{\Wdotar}[1]%
{\hspace{\SOURCE\unitlength}%
\Hcase{\WDOTAR}{#1}{\ARROWLENGTH}}%

\newcommand{\wdotaR}[1]%
{\hspace{\SOURCE\unitlength}%
\hcasE{\WDOTAR}{#1}{\ARROWLENGTH}}%

\newcommand{\wmono}%
{\hspace{\SOURCE\unitlength}%
\hcase{\WMONO}{\ARROWLENGTH}}%

\newcommand{\Wmono}[1]%
{\hspace{\SOURCE\unitlength}%
\Hcase{\WMONO}{#1}{\ARROWLENGTH}}%

\newcommand{\wmonO}[1]%
{\hspace{\SOURCE\unitlength}%
\hcasE{\WMONO}{#1}{\ARROWLENGTH}}%

\newcommand{\wepi}%
{\hspace{\SOURCE\unitlength}%
\hcase{\WEPI}{\ARROWLENGTH}}%

\newcommand{\Wepi}[1]%
{\hspace{\SOURCE\unitlength}%
\Hcase{\WEPI}{#1}{\ARROWLENGTH}}%

\newcommand{\wepI}[1]%
{\hspace{\SOURCE\unitlength}%
\hcasE{\WEPI}{#1}{\ARROWLENGTH}}%

\newcommand{\wbimo}%
{\hspace{\SOURCE\unitlength}%
\hcase{\WBIMO}{\ARROWLENGTH}}%

\newcommand{\Wbimo}[1]%
{\hspace{\SOURCE\unitlength}%
\Hcase{\WBIMO}{#1}{\ARROWLENGTH}}%

\newcommand{\wbimO}[1]%
{\hspace{\SOURCE\unitlength}%
\hcasE{\WBIMO}{#1}{\ARROWLENGTH}}%

\newcommand{\wiso}%
{\hspace{\SOURCE\unitlength}%
\Hcase{\WAR}{\cong}{\ARROWLENGTH}}%

\newcommand{\Wiso}[1]%
{\hspace{\SOURCE\unitlength}%
\Hcase{\WAR}{#1}{\ARROWLENGTH}}%

\newcommand{\wisO}[1]%
{\hspace{\SOURCE\unitlength}%
\hcasE{\WAR}{#1}{\ARROWLENGTH}}%

\newcommand{\wbiar}%
{\hspace{\SOURCE\unitlength}%
\hbicase{\WBIAR}{\ARROWLENGTH}}%

\newcommand{\Wbiar}[2]%
{\hspace{\SOURCE\unitlength}%
\Hbicase{\WBIAR}{#1}{#2}{\ARROWLENGTH}}%

\newcommand{\weql}%
{\hspace{\SOURCE\unitlength}%
\hbicase{\EEQL}{\ARROWLENGTH}}%

\newcommand{\wadjar}%
{\hspace{\SOURCE\unitlength}%
\hbicase{\WADJAR}{\ARROWLENGTH}}%

\newcommand{\Wadjar}[2]%
{\hspace{\SOURCE\unitlength}%
\Hbicase{\WADJAR}{#1}{#2}{\ARROWLENGTH}}%

\newcommand{\Vcase}[3]{\makebox[0pt]%
{\makebox[0pt][r]{\raisebox{0pt}[0pt][0pt]{${#2}\hspace{2pt}$}}}#1{#3}}%

\newcommand{\Vbicase}[4]{\makebox[0pt]%
{\makebox[0pt][r]{\raisebox{0pt}[0pt][0pt]{$#2$\hspace{4pt}}}#1{#4}%
\makebox[0pt][l]{\raisebox{0pt}[0pt][0pt]{\hspace{5pt}$#3$}}}}%

\newcommand{\SAR}[1]%
{\begin{picture}(0,0)%
\put(0,0){\makebox(0,0)%
{\begin{picture}(0,#1)%
\put(0,#1){\vector(0,-1){#1}}%
\end{picture}}}\end{picture}}%

\newcommand{\SDOTAR}[1]%
{\truex{100}\truey{300}%
\NUMBEROFDOTS=#1%
\divide\NUMBEROFDOTS by \value{y}%
\begin{picture}(0,0)%
\put(0,0){\makebox(0,0)%
{\begin{picture}(0,#1)%
\multiput(0,#1)(0,-\value{y}){\NUMBEROFDOTS}%
{\circle*{\value{x}}}%
\put(0,0){\vector(0,-1){0}}%
\end{picture}}}\end{picture}}%

\newcommand{\SMONO}[1]%
{\begin{picture}(0,0)%
\put(0,0){\makebox(0,0)%
{\begin{picture}(0,#1)%
\put(0,#1){\vector(0,-1){#1}}%
\truex{300}\truey{600}%
\put(0,#1){\begin{picture}(0,0)%
\put(-\value{x},-\value{x}){\line(1,0){\value{y}}}\end{picture}}%
\end{picture}}}\end{picture}}%

\newcommand{\SEPI}[1]%
{\begin{picture}(0,0)%
\put(0,0){\makebox(0,0)%
{\begin{picture}(0,#1)%
\put(0,#1){\vector(0,-1){#1}}%
\truex{300}\truey{600}\truez{800}%
\put(-\value{x},\value{z}){\line(1,0){\value{y}}}%
\end{picture}}}\end{picture}}%

\newcommand{\SBIMO}[1]%
{\begin{picture}(0,0)%
\put(0,0){\makebox(0,0)%
{\begin{picture}(0,#1)%
\put(0,#1){\vector(0,-1){#1}}%
\truex{300}\truey{600}\truez{800}%
\put(0,#1){\begin{picture}(0,0)%
\put(-\value{x},-\value{x}){\line(1,0){\value{y}}}\end{picture}}%
\put(-\value{x},\value{z}){\line(1,0){\value{y}}}%
\end{picture}}}\end{picture}}%

\newcommand{\SBIAR}[1]%
{\begin{picture}(0,0)%
\truex{200}%
\put(0,0){\makebox(0,0)%
{\begin{picture}(0,#1)\put(-\value{x},#1){\vector(0,-1){#1}}%
\put(\value{x},#1){\vector(0,-1){#1}}%
\end{picture}}}\end{picture}}%

\newcommand{\SEQL}[1]%
{\begin{picture}(0,0)%
\truex{100}%
\put(0,0){\makebox(0,0)%
{\begin{picture}(0,#1)\put(-\value{x},#1){\line(0,-1){#1}}%
\put(\value{x},#1){\line(0,-1){#1}}%
\end{picture}}}\end{picture}}%

\newcommand{\Sarv}[2]{\Vcase{\SAR}{#1}{#200}}%

\newcommand{\Sar}[1]{\Sarv{#1}{50}}%

\newcommand{\Sisov}[2]%
{\Vbicase{\SAR}{#1\hspace{-2pt}}{\hspace{-2pt}\cong}{#200}}%

\newcommand{\NAR}[1]%
{\begin{picture}(0,0)%
\put(0,0){\makebox(0,0)%
{\begin{picture}(0,#1)\put(0,0){\vector(0,1){#1}}%
\end{picture}}}\end{picture}}%

\newcommand{\NDOTAR}[1]%
{\truex{100}\truey{300}%
\NUMBEROFDOTS=#1%
\divide\NUMBEROFDOTS by \value{y}%
\begin{picture}(0,0)%
\put(0,0){\makebox(0,0)%
{\begin{picture}(0,#1)%
\multiput(0,0)(0,\value{y}){\NUMBEROFDOTS}%
{\circle*{\value{x}}}%
\put(0,#1){\vector(0,1){0}}%
\end{picture}}}\end{picture}}%

\newcommand{\NMONO}[1]%
{\begin{picture}(0,0)%
\put(0,0){\makebox(0,0)%
{\begin{picture}(0,#1)%
\put(0,0){\vector(0,1){#1}}%
\truex{300}\truey{600}%
\put(-\value{x},\value{x}){\line(1,0){\value{y}}}%
\end{picture}}}%
\end{picture}}%

\newcommand{\NEPI}[1]%
{\begin{picture}(0,0)%
\put(0,0){\makebox(0,0)%
{\begin{picture}(0,#1)%
\put(0,0){\vector(0,1){#1}}%
\truex{300}\truey{600}\truez{800}%
\put(0,#1){\begin{picture}(0,0)%
\put(-\value{x},-\value{z}){\line(1,0){\value{y}}}\end{picture}}%
\end{picture}}}\end{picture}}%

\newcommand{\NBIMO}[1]%
{\begin{picture}(0,0)%
\put(0,0){\makebox(0,0)%
{\begin{picture}(0,#1)%
\put(0,0){\vector(0,1){#1}}%
\truex{300}\truey{600}\truez{800}%
\put(-\value{x},\value{x}){\line(1,0){\value{y}}}%
\put(0,#1){\begin{picture}(0,0)%
\put(-\value{x},-\value{z}){\line(1,0){\value{y}}}\end{picture}}%
\end{picture}}}\end{picture}}%

\newcommand{\NBIAR}[1]%
{\begin{picture}(0,0)%
\truex{200}%
\put(0,0){\makebox(0,0)%
{\begin{picture}(0,#1)\put(-\value{x},0){\vector(0,1){#1}}%
\put(\value{x},0){\vector(0,1){#1}}%
\end{picture}}}\end{picture}}%

\newcommand{\Nisov}[2]%
{\Vbicase{\NAR}{#1\hspace{-2pt}}{\hspace{-2pt}\cong}{#200}}%

\newcommand{\NEDOTAR}%
{\truex{100}\truey{212}%
\NUMBEROFDOTS=5800%
\divide\NUMBEROFDOTS by \value{y}%
\begin{picture}(0,0)%
\multiput(-2900,-2900)(\value{y},\value{y}){\NUMBEROFDOTS}%
{\circle*{\value{x}}}%
\put(2900,2900){\vector(1,1){0}}%
\end{picture}}%

\newcommand{\SWDOTAR}%
{\truex{100}\truey{212}%
\NUMBEROFDOTS=5800%
\divide\NUMBEROFDOTS by \value{y}%
\begin{picture}(0,0)%
\multiput(2900,2900)(-\value{y},-\value{y}){\NUMBEROFDOTS}%
{\circle*{\value{x}}}%
\put(-2900,-2900){\vector(-1,-1){0}}%
\end{picture}}%

\newcommand{\SEDOTAR}%
{\truex{100}\truey{212}%
\NUMBEROFDOTS=5800%
\divide\NUMBEROFDOTS by \value{y}%
\begin{picture}(0,0)%
\multiput(-2900,2900)(\value{y},-\value{y}){\NUMBEROFDOTS}%
{\circle*{\value{x}}}%
\put(2900,-2900){\vector(1,-1){0}}%
\end{picture}}%

\newcommand{\NWDOTAR}%
{\truex{100}\truey{212}%
\NUMBEROFDOTS=5800%
\divide\NUMBEROFDOTS by \value{y}%
\begin{picture}(0,0)%
\multiput(2900,-2900)(-\value{y},\value{y}){\NUMBEROFDOTS}%
{\circle*{\value{x}}}%
\put(-2900,2900){\vector(-1,1){0}}%
\end{picture}}%

\newcommand{\ENEAR}[2]%
{\makebox[0pt]{\begin{picture}(0,0)%
\put(0,-150){\makebox(0,0){\begin{picture}(0,0)%
\put(-6600,-3300){\vector(2,1){13200}}%
\truex{200}\truey{800}\truez{600}%
\put(-\value{x},\value{x}){\makebox(0,\value{z})[r]{${#1}$}}%
\put(\value{x},-\value{y}){\makebox(0,\value{z})[l]{${#2}$}}%
\end{picture}}}\end{picture}}}%

\newcommand{\ESEAR}[2]%
{\makebox[0pt]{\begin{picture}(0,0)%
\put(0,-150){\makebox(0,0){\begin{picture}(0,0)%
\put(-6600,3300){\vector(2,-1){13200}}%
\truex{200}\truey{800}\truez{600}%
\put(\value{x},\value{x}){\makebox(0,\value{z})[l]{${#1}$}}%
\put(-\value{x},-\value{y}){\makebox(0,\value{z})[r]{${#2}$}}%
\end{picture}}}\end{picture}}}%

\newcommand{\WNWAR}[2]%
{\makebox[0pt]{\begin{picture}(0,0)%
\put(0,-150){\makebox(0,0){\begin{picture}(0,0)%
\put(6600,-3300){\vector(-2,1){13200}}%
\truex{200}\truey{800}\truez{600}%
\put(\value{x},\value{x}){\makebox(0,\value{z})[l]{${#1}$}}%
\put(-\value{x},-\value{y}){\makebox(0,\value{z})[r]{${#2}$}}%
\end{picture}}}\end{picture}}}%

\newcommand{\WSWAR}[2]%
{\makebox[0pt]{\begin{picture}(0,0)%
\put(0,-150){\makebox(0,0){\begin{picture}(0,0)%
\put(6600,3300){\vector(-2,-1){13200}}%
\truex{200}\truey{800}\truez{600}%
\put(-\value{x},\value{x}){\makebox(0,\value{z})[r]{${#1}$}}%
\put(\value{x},-\value{y}){\makebox(0,\value{z})[l]{${#2}$}}%
\end{picture}}}\end{picture}}}%

\newcommand{\NNEAR}[2]%
{\raisebox{-1pt}[0pt][0pt]{\begin{picture}(0,0)%
\put(0,0){\makebox(0,0){\begin{picture}(0,0)%
\put(-3300,-6600){\vector(1,2){6600}}%
\truex{100}\truez{600}%
\put(-\value{x},\value{x}){\makebox(0,\value{z})[r]{${#1}$}}%
\put(\value{x},-\value{z}){\makebox(0,\value{z})[l]{${#2}$}}%
\end{picture}}}\end{picture}}}%

\newcommand{\SSWAR}[2]%
{\raisebox{-1pt}[0pt][0pt]{\begin{picture}(0,0)%
\put(0,0){\makebox(0,0){\begin{picture}(0,0)%
\put(3300,6600){\vector(-1,-2){6600}}%
\truex{100}\truez{600}%
\put(-\value{x},\value{x}){\makebox(0,\value{z})[r]{${#1}$}}%
\put(\value{x},-\value{z}){\makebox(0,\value{z})[l]{${#2}$}}%
\end{picture}}}\end{picture}}}%

\newcommand{\SSEAR}[2]%
{\raisebox{-1pt}[0pt][0pt]{\begin{picture}(0,0)%
\put(0,0){\makebox(0,0){\begin{picture}(0,0)%
\put(-3300,6600){\vector(1,-2){6600}}%
\truex{200}\truez{600}%
\put(\value{x},\value{x}){\makebox(0,\value{z})[l]{${#1}$}}%
\put(-\value{x},-\value{z}){\makebox(0,\value{z})[r]{${#2}$}}%
\end{picture}}}\end{picture}}}%

\newcommand{\NNWAR}[2]%
{\raisebox{-1pt}[0pt][0pt]{\begin{picture}(0,0)%
\put(0,0){\makebox(0,0){\begin{picture}(0,0)%
\put(3300,-6600){\vector(-1,2){6600}}%
\truex{200}\truez{600}%
\put(\value{x},\value{x}){\makebox(0,\value{z})[l]{${#1}$}}%
\put(-\value{x},-\value{z}){\makebox(0,\value{z})[r]{${#2}$}}%
\end{picture}}}\end{picture}}}%

\newcommand{\Necurve}[2]%
{\begin{picture}(0,0)%
\truex{1300}\truey{2000}\truez{200}%
\put(0,\value{x}){\oval(#200,\value{y})[t]}%
\put(0,\value{x}){\makebox(0,0){\begin{picture}(#200,0)%
\put(#200,0){\vector(0,-1){\value{z}}}%
\put(0,0){\line(0,-1){\value{z}}}\end{picture}}}%
\truex{2500}%
\put(0,\value{x}){\makebox(0,0)[b]{${#1}$}}%
\end{picture}}%

\newcommand{\Nwcurve}[2]%
{\begin{picture}(0,0)%
\truex{1300}\truey{2000}\truez{200}%
\put(0,\value{x}){\oval(#200,\value{y})[t]}%
\put(0,\value{x}){\makebox(0,0){\begin{picture}(#200,0)%
\put(#200,0){\line(0,-1){\value{z}}}%
\put(0,0){\vector(0,-1){\value{z}}}\end{picture}}}%
\truex{2500}%
\put(0,\value{x}){\makebox(0,0)[b]{${#1}$}}%
\end{picture}}%

\newcommand{\Securve}[2]%
{\begin{picture}(0,0)%
\truex{1300}\truey{2000}\truez{200}%
\put(0,-\value{x}){\oval(#200,\value{y})[b]}%
\put(0,-\value{x}){\makebox(0,0){\begin{picture}(#200,0)%
\put(#200,0){\vector(0,1){\value{z}}}%
\put(0,0){\line(0,1){\value{z}}}\end{picture}}}%
\truex{2500}%
\put(0,-\value{x}){\makebox(0,0)[t]{${#1}$}}%
\end{picture}}%

\newcommand{\Swcurve}[2]%
{\begin{picture}(0,0)%
\truex{1300}\truey{2000}\truez{200}%
\put(0,-\value{x}){\oval(#200,\value{y})[b]}%
\put(0,-\value{x}){\makebox(0,0){\begin{picture}(#200,0)%
\put(#200,0){\line(0,1){\value{z}}}%
\put(0,0){\vector(0,1){\value{z}}}\end{picture}}}%
\truex{2500}%
\put(0,-\value{x}){\makebox(0,0)[t]{${#1}$}}%
\end{picture}}%

\newcommand{\Escurve}[2]%
{\begin{picture}(0,0)%
\truex{1400}\truey{2000}\truez{200}%
\put(\value{x},0){\oval(\value{y},#200)[r]}%
\put(\value{x},0){\makebox(0,0){\begin{picture}(0,#200)%
\put(0,0){\vector(-1,0){\value{z}}}%
\put(0,#200){\line(-1,0){\value{z}}}\end{picture}}}%
\truex{2500}%
\put(\value{x},0){\makebox(0,0)[l]{${#1}$}}%
\end{picture}}%

\newcommand{\Encurve}[2]%
{\begin{picture}(0,0)%
\truex{1400}\truey{2000}\truez{200}%
\put(\value{x},0){\oval(\value{y},#200)[r]}%
\put(\value{x},0){\makebox(0,0){\begin{picture}(0,#200)%
\put(0,0){\line(-1,0){\value{z}}}%
\put(0,#200){\vector(-1,0){\value{z}}}\end{picture}}}%
\truex{2500}%
\put(\value{x},0){\makebox(0,0)[l]{${#1}$}}%
\end{picture}}%

\newcommand{\Wscurve}[2]%
{\begin{picture}(0,0)%
\truex{1300}\truey{2000}\truez{200}%
\put(-\value{x},0){\oval(\value{y},#200)[l]}%
\put(-\value{x},0){\makebox(0,0){\begin{picture}(0,#200)%
\put(0,0){\vector(1,0){\value{z}}}%
\put(0,#200){\line(1,0){\value{z}}}\end{picture}}}%
\truex{2400}%
\put(-\value{x},0){\makebox(0,0)[r]{${#1}$}}%
\end{picture}}%

\newcommand{\Wncurve}[2]%
{\begin{picture}(0,0)%
\truex{1300}\truey{2000}\truez{200}%
\put(-\value{x},0){\oval(\value{y},#200)[l]}%
\put(-\value{x},0){\makebox(0,0){\begin{picture}(0,#200)%
\put(0,0){\line(1,0){\value{z}}}%
\put(0,#200){\vector(1,0){\value{z}}}\end{picture}}}%
\truex{2400}%
\put(-\value{x},0){\makebox(0,0)[r]{${#1}$}}%
\end{picture}}%

\newcount\SCALE%

\newcount\NUMBER%

\newcount\LINE%

\newcount\COLUMN%

\newcount\WIDTH%

\newcount\SOURCE%

\newcount\ARROW%

\newcount\TARGET%

\newcount\ARROWLENGTH%

\newcount\NUMBEROFDOTS%

\newcounter{x}%

\newcounter{y}%

\newcounter{z}%

\newcounter{horizontal}%

\newcounter{vertical}%

\newskip\itemlength%

\newskip\firstitem%

\newskip\seconditem%

\newcommand{\printarrow}{}%

\newcommand{\truex}[1]{%
\NUMBER=#1%
\multiply\NUMBER by 100%
\divide\NUMBER by \SCALE%
\setcounter{x}{\NUMBER}}%

\newcommand{\truey}[1]{%
\NUMBER=#1%
\multiply\NUMBER by 100%
\divide\NUMBER by \SCALE%
\setcounter{y}{\NUMBER}}%

\newcommand{\truez}[1]{%
\NUMBER=#1%
\multiply\NUMBER by 100%
\divide\NUMBER by \SCALE%
\setcounter{z}{\NUMBER}}%

\newcommand{\changecounters}[1]{%
\SOURCE=\ARROW%
\ARROW=\TARGET%
\settowidth{\itemlength}{#1}%
\ifdim \itemlength > 2800\unitlength%
\addtolength{\itemlength}{-2800\unitlength}%
\TARGET=\itemlength%
\divide\TARGET by 1310%
\multiply\TARGET by 100%
\divide\TARGET by \SCALE%
\else%
\TARGET=0%
\fi%
\ARROWLENGTH=5000%
\advance\ARROWLENGTH by -\SOURCE%
\advance\ARROWLENGTH by -\TARGET%
\advance\SOURCE by -\TARGET}%

\newcommand{\initialize}[1]{%
\LINE=0%
\COLUMN=0%
\WIDTH=0%
\ARROW=0%
\TARGET=0%
\changecounters{#1}%
\renewcommand{\printarrow}{#1}%
\begin{center}%
\vspace{10pt}%
\begin{picture}(0,0)}%

\newcommand{\DIAGV}[2]{%
\SCALE=#1%
\setlength{\unitlength}{655sp}%
\multiply\unitlength by \SCALE%
\divide\unitlength by 100%
\initialize{\mbox{$#2$}}}%

\newcommand{\n}[1]{%
\changecounters{\mbox{$#1$}}%
\put(\COLUMN,\LINE){\makebox(0,0){\printarrow}}%
\thinlines%
\renewcommand{\printarrow}{\mbox{$#1$}}%
\advance\COLUMN by 4000}%

\newcommand{\nn}[1]{%
\put(\COLUMN,\LINE){\makebox(0,0){\printarrow}}%
\thinlines%
\ifnum \WIDTH < \COLUMN%
\WIDTH=\COLUMN%
\else%
\fi%
\advance\LINE by -4000%
\COLUMN=0%
\ARROW=0%
\TARGET=0%
\changecounters{\mbox{$#1$}}%
\renewcommand{\printarrow}{\mbox{$#1$}}}%

\newcommand{\conclude}{%
\put(\COLUMN,\LINE){\makebox(0,0){\printarrow}}%
\thinlines%
\ifnum \WIDTH < \COLUMN%
\WIDTH=\COLUMN%
\else%
\fi%
\setcounter{horizontal}{\WIDTH}%
\setcounter{vertical}{-\LINE}%
\end{picture}}%

\newcommand{\diag}{%
\conclude%
\raisebox{0pt}[0pt][\value{vertical}\unitlength]{}%
\hspace*{\value{horizontal}\unitlength}%
\vspace{10pt}%
\end{center}%
\setlength{\unitlength}{1pt}}%

\newcommand{\diagv}[3]{%
\conclude%
\NUMBER=#1%
\rule{0pt}{\NUMBER pt}%
\hspace*{-#2pt}%
\raisebox{0pt}[0pt][\value{vertical}\unitlength]{}%
\hspace*{\value{horizontal}\unitlength}
\NUMBER=#3%
\advance\NUMBER by 10%
\vspace*{\NUMBER pt}%
\end{center}%
\setlength{\unitlength}{1pt}}%

\newcommand{\N}[1]%
{\raisebox{0pt}[7pt][0pt]{$#1$}}%

\newcommand{\crosslength}[2]{%
\settowidth{\firstitem}{#1}%
\settowidth{\seconditem}{#2}%
\ifdim\firstitem < \seconditem%
\itemlength=\seconditem%
\else%
\itemlength=\firstitem%
\fi%
\divide\itemlength by 2%
\hspace{\itemlength}}%

\newtheorem{proposition}{Proposition}

\newtheorem{definition}{Definition}

\def\N{I\!\!N}

\def\Z{Z\!\!\!Z}

\newcommand{\cue}{I\!\! F}

\newcommand{\Na}{I\!\! N}
\def\vae{\varepsilon}

\begin{document}
\title{GROUP CONVOLUTIONAL CODES}
\author{
S. Estrada \\Dept. de Matem\'atica Aplicada.\\ Universidad de
Murcia. Campus de Espinardo \\ Espinardo, Murcia,  30100, Spain \\
 \and J. R. Garc\'{\i}a Rozas, \and J. Peralta, \and
E. S\'anchez Garc\'{\i}a \and   Dept. \'{A}lgebra y An\'{a}lisis
Matem\'{a}tico\\ Universidad de Almer\'{\i}a \\ 04120
Almer\'{\i}a, Spain \\ e-mail: jrgrozas@ual.es }

\date{}
\maketitle
\begin{abstract}

In this note we introduce the concept of group convolutional code.
We make a complete classification of the minimal $S_3$-convolutional
codes over the field of five elements by means of Jategaonkar's
theorems.

\bigskip

\medskip

\end{abstract}

\section{Introduction}
Block codes as left ideals in group algebras were introduced by S.
D. Bermann in \cite{berman}. After that, several papers of
MacWilliams, Landrock, Damgard, Lieber, Ward, Zimmermman and
others gave more credit to this theory (\cite{dam1}, \cite{land1},
\cite{lieb1},\cite{mac1},\cite{PERA},\cite{zimmer2}). In the
context of convolutional codes, P. Piret \cite{PIR2}, studied the
$H$-codes, which can be seen as a generalized version of the group
block codes in the convolutional case.

On the other hand, the concept of cyclic convolutional codes and
their first properties were proposed by P. Piret and C. Roos in
\cite{PIR1} and \cite{ROOS}, respectively. More recently, H.
Gluesing-Luerssen et al. (\cite{OCCCL}, \cite{OCCc}) continue the
study of cyclic convolutional codes. In the present paper, we give
a definition of group convolutional code, which is a
generalization of cyclic convolutional code and group  block code.
We introduce some important techniques in non-commutative algebra,
concretely, the structure theorems of skew polynomials rings given
in \cite{JATE} by Jategaonkar.

The paper is organized as follows. In Section 2 we make the
necessaries definitions related with convolutional codes we will use
throughout the paper. Then we introduce the concept of group
convolutional code and minimal one, this last will be the main
object of our study since they are the building blocks for the rest
of the codes. Next we summarize Jategaonkar's result on the
structure of skew polynomial rings over semisimple rings, that we
will use in the last section. Finally, Section 3 deals with the
classification of the minimal $S_3$-convolutional codes over the
field of five elements. The isomorphism established between the skew
polynomial ring and certain direct sums of rings of matrices over
simplest skew polynomial rings will be crucial. Note that these
codes are the smallest non-commutative group convolutional codes to
consider. This result opens the way to consider more complicated
examples.

\section{Preliminaries and first results}

Throughout this paper, $\cue$ denotes a finite field and $n$ a
positive integer such that the characteristic of $\cue $, $char(\cue
)$, does not divide $n$. This assumption guarantees that for any
group $G$ of order $n$, the group algebra $\cue [G]$ is semi-simple.

This paper deals with convolutional codes with additional
algebraic structure. We adopt the following definition of
convolutional code from \cite{OCCc}.

\begin{definition} A convolutional code of length $n$ and dimension $k$
is a direct summand ${\cal C}$ of $\cue[z]^n$ of rank $k$ as $\cue
[z]$-module.
\end{definition}

Let $r$ be a positive integer. Any matrix $M\in M_{r\times n}(\cue
[z])$ with rows given by a generating set of ${\cal C}$ as $\cue
[z]$-module is called generating matrix of the code ${\cal C}$. If
$r=k$, then $M$ is called generator matrix or encoder of ${\cal
C}$.

The maximal degree of the $k$-minors of an encoder $M$ is called
the {\it complexity} of the code. A code of complexity zero is
said to be a block code.

The {\it free distance} of a convolutional code is defined as
follows. First, given $v=\sum_{i=0}^mv_iz^i\in \cue[z]^n$ where
$v_j\in \cue^n$, we define its weight as
$wt(v)=\sum_{i=0}^mwt(v_i)$, where $wt(v_i)$ is the usual Hamming
weight of the vector $v_i\in \cue^n$. Then, the free distance of a
convolutional code ${\cal C}\subseteq \cue[z]^n$ is defined as,
$dist({\cal C})=min\{ wt(v)\; | \; v\in {\cal C}-\{ 0\} \} $.

We call  $(n,k,\delta )$-convolutional code a code with length
$n$, dimension (or rank) $k$ and complexity $\delta$. We say that
a $(n,k,\delta )$-convolutional code with free distance $d$,
${\cal C}$, is a MDS code (maximal distance separable) if
$d=S(n,k,\delta )$, where $S(n,k,\delta )$ is the generalized
Singleton bound, $S(n,k,\delta) =(n-k)(\lfloor
\frac{\delta}{k}\rfloor +1 )+\delta +1$. For a given size field
$q$, we have the so called {\it Griesmer bound} for convolutional
codes over the field of $q$ elements. It is defined as

\noindent $G(n,k, \delta ; m)_q=  max \left \{ d' \in \{ 1,...,
S(n,k,\delta )\} \; \mid \; \sum_{l=0}^{k(m+i)-\delta -1} \lceil
\frac{d'}{q^l}\rceil \leq n(m+i)\right .$ $\left . \mbox{for all}
\; i\in \widehat{\Na}\right \} .$

 Here $m$ denotes the maximum taken over
the Forney indices of a $(n,k,\delta )$-convolutional code, and it
is called the {\it memory} of the code. Also, $\widehat{\Na}$
denotes $\{ 1,2,...\}$ if $km=\delta $ or $\{ 0,1,2,...\}$ if $km
> \delta$. A convolutional code over a field of $q$ elements is
said to be optimal if it reaches the Griesmer bound (see
\cite{DCCC} ).
\medskip

Let $G=\{ g_1,...,g_n\}$ be a finite group of order $n$. We consider
the group $\cue$-algebra $A=\cue [G]$ and the $\cue$-isomorphism
$\beta :\cue^n\rightarrow A$ given by $\beta
(v_1,...,v_n):=\sum_{i=1}^nv_ig_i$. On the other hand, we have the
canonical isomorphism $\psi:\cue [z]^n\rightarrow \cue^n[z]$. Given
$y\in\cue [z]^n$, let $\psi (y)=\sum_{j\geq 0} z^jw_j\in\cue^n[z]$.
Then we define $\rho :\cue[z]^n\rightarrow A[z]$ by $\rho
(y)=\sum_{j\geq 0} z^j\beta(w_j).$ It is clear that $\rho$ is a
$\cue [z]$-isomorphism. We identify the $\cue [z]$-submodules of
$\cue [z]^n$ with the $\cue [z]$-submodules of $A[z]$ via $\rho$.

Now, let $\sigma$ be an $\cue$-automorphism of $A$ and ${\cal
R}=A[z;\sigma ]$ be the skew polynomial ring. The multiplication
rule in ${\cal R}$ is given by $az=z\sigma (a)$ for all $a\in A$.
The map $\rho_{\sigma}: \cue[z]^n\rightarrow A[z;\sigma ]$ defined
just like $\rho$ is the key for the next definitions ( in
\cite{OCCCL} essentially appears the respective definitions in the
particular case of a cyclic group). Note that $\rho_{\sigma}$ is
an isomorphism of left $\cue [z]$-modules.

\begin{definition} Let ${\cal C}\subseteq \cue[z]^n$ be a convolutional code.
We say that ${\cal C}$ is a $(G,\sigma )$-convolutional code if
$\rho_{\sigma}({\cal C})$ is a direct summand left ideal of ${\cal
R}$.
\end{definition}

We will see that this definition coincides with the usual one
where only is required that $\rho_{\sigma}({\cal C})$ is a direct
summand as $\cue [z]$-module.

\begin{proposition} Let ${\cal C}\subseteq \cue[z]^n$ be a convolutional code.
The following conditions are equivalent.

a) ${\cal C}$ is a $(G,\sigma )$-convolutional code.

b) $\rho_{\sigma}({\cal C})$ is a left ideal of ${\cal R}$ and there
is an $\cue [z]$-submodule $K$ of ${\cal R}$ such that
$\rho_{\sigma}({\cal C})\oplus K= {\cal R}$.
\end{proposition}
{\bf Proof. } $a)\Rightarrow b)$ is obvious since any left ideal
of ${\cal R}$ is, in particular, a $\cue [z]$-submodule.

$b)\Rightarrow a)$ Suppose $\rho_{\sigma}({\cal C})\oplus K = {\cal
R}$ as $\cue [z]$-modules. Then there is an $\cue [z]$-linear map
$\pi : {\cal R}\rightarrow \rho_{\sigma}({\cal C})$ such that $\pi
(x)=x$ for all $x\in \rho_{\sigma}({\cal C})$. Define $\overline
{\pi}: {\cal R}\rightarrow \rho_{\sigma}({\cal C})$ by
$\overline{\pi}(a)=\frac{1}{m}(\sum_{d\in U}d\ \pi (d^{-1}a))$,
where $U=U(\cue [G])$ is the group of units of $\cue [G]$ and $m$ is
its order. It is clear that $\overline{\pi}(x)=x$ for all $x\in
\rho_{\sigma}({\cal C})$. We will show that $\overline{\pi}$ is
${\cal R}$-linear and so $\rho_{\sigma}({\cal C})$ would be a direct
summand of ${\cal R}$ as left ${\cal R}$-modules. It is enough to
prove that $\overline{\pi}(ha)=h\overline{a}$ and
$\overline{\pi}(za)=z\overline{\pi}(a)$ for all $h\in G$, $a\in
{\cal R}$. Now, $\overline{\pi}(ha)=\frac{1}{m}(\sum_{d\in U}d\ \pi
(d^{-1}ha))=\frac{1}{m}(\sum_{d\in U}hh^{-1}d\ \pi
(d^{-1}ha))=h(\frac{1}{m}(\sum_{d\in U}h^{-1}d\ \pi
(d^{-1}ha)))=h\overline{\pi}(a).$

Also, $\overline{\pi}(za)=\frac{1}{m}(\sum_{d\in U}d\ \pi
(d^{-1}za))=z(\frac{1}{m}(\sum_{d\in U}\sigma (d)\ \pi (\sigma
(d)^{-1}a)))=z\overline{\pi}(a)$,  (the last equality holds
because $\sigma$ produces a permutation on the elements in $U$).
$\Box$

\begin{definition}
We say that a  $(G,\sigma )$-convolutional code ${\cal C}$ is
minimal if $\rho_{\sigma}({\cal C})$ is indecomposable as left
$A[z;\sigma ]$-module.
\end{definition}

\begin{proposition}
a) Any minimal $(G,\sigma )$-convolutional  code does not contain
any other proper $(G,\sigma )$-convolutional code.

b) Any $(G,\sigma )$-convolutional code is a direct sum of minimal
$(G,\sigma )$-con\-vo\-lu\-tio\-nal codes.
\end{proposition}
{\bf Proof.} a) Let ${\cal C}$ be a minimal $(G,\sigma
)$-convolutional  code and ${\cal L}\subseteq {\cal C}$ a $(G,\sigma
)$-convolutional  code different from ${\cal C}$. Then ${\cal
R}=\rho_{\sigma}({\cal L})\oplus K$ as left ${\cal R}$-modules for
some $K\leq {\cal R}$. This implies that $\rho_{\sigma}({\cal
C})=\rho_{\sigma}({\cal L})\oplus (K\cap \rho_{\sigma}({\cal C}))$
which is a contradiction with the minimality of ${\cal C}$.

b) Let ${\cal C}$ be a $(G,\sigma )$-convolutional code. Then
$I=\rho_{\sigma}({\cal C})$ is a direct summand left ideal of ${\cal
R}$. If $I$ is indecomposable then it is done. In the contrary case,
$I=I_1\oplus I_2$ where $I_i$ is a nonzero left ideal of ${\cal R}$
for $i=1,2$. Again if both $I_i$ are indecomposable it is done. This
procedure can be repeated and must stop since the ideal $I$ has
finite rank as $\cue [z]$-module and the $I_i$'s are free $\cue
[z]$-modules.

\bigskip

It is standard that any minimal $(G,\sigma )$-convolutional code is
generated as left $A[z;\sigma ]$-module by a primitive idempotent
element of $A[z;\sigma ]$. This paper mainly deals with the problem
of finding these primitive idempotents. We are interested in the
matrix approach of $A[z;\sigma ]$. Next, we make an account of
results on the interpretation of the elements of $A[z;\sigma ]$ as
matrices in some matrix ring. We use Jategaonkar's results (cf.
\cite{JATE}) in order to give an explicit isomorphism of rings
between $A[z;\sigma ]$ and the rings constructed via matrix rings.

For the rest of this section, let $A$ be a finite ring (non
necessarily commutative), $\sigma :A\rightarrow A$ be an
automorphism and $z$ an indeterminate. The skew polynomial ring
${\cal R}=A[z;\sigma]$ admits a variable change in $z$ such that
${\cal R}$ is again a skew polynomial ring: let $u$ be a unit in $A$
and $\overline{u}$ the inner automorphism of $A$ defined by
$\overline{u}(a)=u^{-1}au$, $a\in A$. It is easy to check that
$A[z;\sigma ]=A[zu;\overline{u}\sigma ]$.

The following rings are intimately related to the skew polynomial
rings. Let $K$ be a ring and $\rho :K\rightarrow K$ an
automorphism. Let $D=K[x;\rho ]$, $m>0$ and $P$ the subring of
$M_m(D)$ consisting of all the matrices $(d_{ij})$ satisfying the
next two conditions: (1) $d_{ij}\in D$ $\forall i,j$; (2)
$d_{ij}\in xD$ if $i>j$. We denote the subring $P$ by $\{ K,m,\rho
,x\}$. We also denote by $I_n$ the set $\{ 1,...,n\}$.

We recall the concept of {\it set of matrix units} that appears is
\cite[P. 52]{jacob}. Let $A$ be a ring. A finite subset $\{
e_{ij}\, : \, i,j\in I_n\}$ in $A$ is called set of matrix units
in $A$ if verifies the following two conditions:
$$\sum_{i=1}^ne_{ii}=1\;\;\; \mbox{and}\;\;\;
e_{ij}e_{kl}=\delta_{jk}e_{il}$$ where $\delta_{jk}$ is the
Kronecker delta.  In particular, $e_{ij}\neq 0$ for all $i,j\in
I_n$.

A central idempotent element $f$ in $A$ is called {\it
semiprimitive} if $f$ is pri\-mi\-ti\-ve in the center of $A$.

The following fact will be used frequently in the next section. Let
$A$ be a semisimple finite ring and $\{ f_1,...,f_m\}$ a complete
set of semiprimitive idempotent elements in $A$. Assume that $\sigma
:A\rightarrow A$ is an automorphism such that $\sigma (f_i)=f_{\pi
(i)}$ where $\pi$ is the cycle over $I_m$ given by $\pi =(1\;
2...m)$. Let ${\cal R}=A[z;\sigma ]$. Then, by \cite[Lemma
3.1]{JATE}, there exists a finite field $K$, an automorphism $\rho
:K\rightarrow K$ and a positive integer $n$ such that ${\cal R}\cong
M_n(\{ K,m,\rho ,x\} )$ for some indeterminate $x$.

Note that the above positive integer $n$ is the cardinality of a
complete set of matrix units in $Af_1$.

\section{$S_3$-convolutional codes}

In this section we are going to determinate the minimal
$S_3$-convolutional codes over the field with five elements via
Jategaonkar's theorems \cite{JATE}. We fix the field with $5$
elements $\cue_5$ and let $A=\cue_5[S_3]$. The ring $A$ is
semisimple by Maschke Theorem. First, we calculate a complete set
of primitive orthogonal idempotents elements of $A$ by means of
theory of Young diagrams (see \cite[pg. 190]{kup}). The list of
the four idempotent is the following:
\medskip

$\vae_1=I+(1\, 2)+(1\, 3)+(2\, 3)+(1\, 2\, 3)+(1\, 3\, 2),$
\medskip

$\vae_2=I+4(1\, 2)+4(1\, 3)+4(2\, 3)+(1\, 2\, 3)+(1\, 3\, 2),$
\medskip

$\vae_3=2I+3(1\, 2)+2(2\, 3)+3(1\, 2\, 3),$
\medskip

$\vae_4=2I+2(1\, 2)+3(2\, 3)+3(1\, 3\, 2).$

Then $A=\vae_1A\oplus \vae_2A\oplus\vae_3A\oplus \vae_4A$, where
$\vae_1A\cong \vae_2A\cong \cue_5$ and $\vae_3A\oplus\vae_4A\cong
M_2(\cue_5)$ as rings. The corresponding semiprimitive idempotents
are $f_1=\vae_1$, $f_2=\vae_2$ and $f_3=\vae_3+\vae_4$.

We consider two classes of $\cue_5$-automorphism of $A$ attending
to the feasible permutation that produces over the set $\{
f_1,f_2,f_3\}$. One class will be represented by the identity
permutation and the other by the permutation $(1\, 2)$. By
\cite[Theorem 3.3]{JATE}, two automorphisms that produce the same
permutation also produce isomorphic skew polynomial rings.
Moreover, we will prove later that they are {\it isometric}, in
the sense that there is ring isomorphisms between them that
preserve the weight of the elements. So we only take in our study
the identity automorphism (for the identity permutation) and any
automorphism $\sigma\in Aut_{\cue_5}(A)$ such that $\sigma
(f_1)=f_2$, $\sigma (f_2)=f_1$ and $\sigma (f_3)=f_3$ (note that
any automorphism maps $f_1$ to $f_1$ or $f_2$, $f_2$ to $f_2$ or
$f_1$ and $f_3$ to $f_3$).

\subsection{The case of the permutation $(1\ 2)$}

We begin with the second type of automorphism. We take the
automorphism $\sigma$ such that $\sigma (I)=I$, $\sigma (1\,
2)=4(1\, 2)$, $\sigma (1\, 3)=4(1\, 3)$, $\sigma (2\, 3)=4(2\,
3)$, $\sigma (1\, 2\, 3)=(1\, 2\, 3)$, $\sigma (1\, 3\, 2)=(1\,
3\, 2)$. It can be checked that $\sigma$ verifies the above
conditions over $\{ f_1,f_2,f_3\}$.

By \cite[Lemma 3.2]{JATE}, $A[z;\sigma ]=Ag_1[zg_1;\sigma_1]\oplus
Ag_2[zg_2;\sigma_2]$, where $g_1=f_1+f_2$, $g_2=f_3$,
$\sigma_1=\sigma |_{Ag_1}$ and $\sigma_2=\sigma |_{Ag_2}$.

Let $b_1=I$, $b_2=(1\, 2)$, $b_3=(1\, 3)$, $b_4=(2\, 3)$, $b_5=(1\,
2\, 3)$ and $b_6=(1\, 3\, 2)$. Given $h\in A[z; \sigma ]$, we have
$$h=\sum_{i=0}^mz^i(\sum_{j=1}^6a_{ij}b_j)=\sum_{j=1}^6(\sum_{i=0}^mz^ia_{ij})b_j$$
with $a_{ij}\in \cue_5$. Then,
$$h=hg_1+hg_2=\sum_{j=1}^6(\sum_{i=0}^m(zg_1)^ia_{ij})b_jg_1+
\sum_{j=1}^6(\sum_{i=0}^m(zg_2)^ia_{ij})b_jg_2.$$ We study
separately $hg_1$ and $hg_2$.

By \cite[Theorem 2.1]{JATE}, there exists an isomorphism  $\phi_1:
Ag_1[zg_1;\sigma_1]\rightarrow {\cal S}$, where ${\cal S}=\{
\cue_5, 2,\rho ,x\}\subseteq M_2(\cue_5 [x;\rho ])$, $\vae_2A\cong
\cue_5$, $x=(zg_1)^2$ and $\rho
=\sigma_1^2=id_{\vae_2A}:\vae_2A\rightarrow \vae_2A$. Hence, the
ring ${\cal S}$ is simply the subring of $M_2(\cue_5[x])$ given by
${\cal S}=\{ \left ( \begin{array}{cc} p_{11} & p_{12}\\
xp_{21} & p_{22}\\ \end{array}\right ) \; \mid\; p_{ij}\in
\cue_5[x]\}$. To understand $\phi_1(hg_1)$ is enough to calculate
$\phi_1((zg_1))$ and $\phi_1(b_jg_1)$, for all  $j\in I_6$. It is
easy to see that $b_jg_1$ is equal to $2(I+(1\, 2\, 3)+(1\, 3\,
2))$ or $2((1\, 2)+(1\, 3)+(2\, 3))$ for all  $j\in I_6$. Then, by
the proof of \cite[Theorem 2.1]{JATE}, we have
$\phi_1 (zg_1)=\left ( \begin{array}{cc} 0 & 1\\
x & 0\\ \end{array}\right )$, $\phi_1(((1\, 2)+(1\, 3)+(2\,
3))g_1)=\left ( \begin{array}{cc} 3 & 0\\
0 & 2\\ \end{array}\right )$ and $\phi_1((I+(1\, 2\, 3)+(1\,
3\, 2))g_1)=\left ( \begin{array}{cc} 3 & 0\\
0 & 3\\ \end{array}\right )$.

Note that $\phi_1(\vae_1g_1)=\left ( \begin{array}{cc} 1 & 0\\
0 & 0\\ \end{array}\right )$, $\phi_1(\vae_2g_1)=\left ( \begin{array}{cc} 0 & 0\\
0 & 1\\ \end{array}\right )$ and by applying $\phi_1$ to the sum
$((1\, 2)+(1\, 3)+(2\, 3))g_1+(I+(1\, 2\, 3)+(1\, 3\,
2))g_1=\vae_1g_1$
we get precisely $\left ( \begin{array}{cc} 3 & 0\\
0 & 2\\ \end{array}\right ) +\left ( \begin{array}{cc} 3 & 0\\
0 & 3\\ \end{array}\right ) =\left ( \begin{array}{cc} 1 & 0\\
0 & 0\\ \end{array}\right )$.

\bigskip

Now we focus our attention on the direct summand
$Ag_2[zg_2;\sigma_2]$. By \cite[Lemma 3.1]{JATE}, there exists an
isomorphism $\psi :Ag_2[zg_2;\sigma_2]\rightarrow
M_2(\cue_5[zg_2u;\overline{u}\sigma_2])$, where $u$ is a unit in
$Ag_2$. We will make effective this isomorphism.

First we find an isomorphism $\delta :Ag_2\rightarrow
M_2(\cue_5)$. Let $\vae_{33}=\vae_3$, $\vae_{44}=\vae_4$,
$\vae_{34}=(1\, 3)\vae_4$ and $\vae_{43}=(1\, 3)\vae_3$. Then, by
the theory of Young diagrams, the set $\{
\vae_{33},\vae_{34},\vae_{43},\vae_{44}\}$ is a set of matrix
units for $Ag_2$ (see \cite{jacob}). Hence the assignation
$\vae_{33}\mapsto \left ( \begin{array}{cc} 1 & 0 \\ 0 & 0 \\
\end{array} \right )$, $\vae_{34}\mapsto \left ( \begin{array}{cc} 0 & 1 \\ 0 & 0 \\
\end{array} \right )$, $\vae_{43}\mapsto \left ( \begin{array}{cc} 0 & 0 \\ 1 & 0 \\
\end{array} \right )$, $\vae_{44}\mapsto \left ( \begin{array}{cc} 0 & 0 \\ 0 & 1 \\
\end{array} \right )$, will produce the isomorphism $\delta$.
Concretely, given $ag_2\in Ag_2$, we define
$a_{ij}=\sum_{k=3}^4\vae_{ki}ag_2\vae_{jk}$ with $j,i\in \{
3,4\}$. Then $a_{ij}$ belongs to the center of $Ag_2$
(\cite{jacob}), $Cent(Ag_2)\cong \cue_5$, and $\delta
(ag_2)=(a_{ij})$ verifies the above.

Now we need to know how $\sigma_2:Ag_2\rightarrow Ag_2$ is induced
in $M_2(\cue_5)$, i.e., we must find an automorphism
$\widehat{\sigma}_2 :M_2(\cue_5)\rightarrow M_2(\cue_5)$ such that
the diagram \DIAGV{80} {Ag_2} \n{\Ear{\sigma_2}}  \n{Ag_2} \nn
{\Sar{\delta}} \n{} \n{\Sar{\delta}} \nn {M_2(\cue_5)}
\n{\Ear{\widehat{\sigma}_2}} \n{M_2(\cue_5)} \diag is commutative.
It is clear that $\widehat{\sigma}_2=\delta \sigma_2\delta^{-1}$.
Then, easy calculations show that
$\widehat{\sigma}_2\left ( \begin{array}{cc} a & b \\ c & d \\
\end{array} \right ) =B\left ( \begin{array}{cc} a & b \\ c & d \\
\end{array} \right ) B^{-1}$, where $B=\left ( \begin{array}{cc} 4 & 3 \\ 2 & 1 \\
\end{array} \right )$.

Hence we have the isomorphism induced by $\delta$ in the obvious
manner:  $\overline{\delta }: Ag_2[zg_2;\sigma_2]\rightarrow
M_2(\cue_5)[y;\widehat{\sigma}_2]$, $zg_2\mapsto y$,
$\overline{\delta}|_{Ag_2}=\delta$. On the other hand,
$M_2(\cue_5)[y;\widehat{\sigma}_2]=M_2(\cue_5)[yB;\overline{B}\widehat{\sigma}_2]=M_2(\cue_5)[yB]$.
Taking $yB=x$, we finally get
$$M_2(\cue_5)[y;\widehat{\sigma}_2]=M_2(\cue_5)[x]\cong
M_2(\cue_5[x]),$$ with the last isomorphism the canonical one. Let
$$\phi_2:Ag_2[zg_2;\sigma_2]\rightarrow M_2(\cue_5[x])$$ be the
composition of $\overline{\delta}$ with the canonical isomorphism.
Then, it is clear that  $\phi_2(zg_2)=xB^{-1}$ and given $\alpha
=a_{33}\vae_{33}+a_{44}\vae_{44}+a_{34}\vae_{34}+a_{43}\vae_{43}\in
Ag_2$,
$\phi_2(\alpha )=\left ( \begin{array}{cc} a_{33} & a_{34} \\ a_{43} & a_{44} \\
\end{array} \right )$ .
\bigskip

Once we have completely described the isomorphisms $\phi_1$ and
$\phi_2$, we have the ring  isomorphism
$\phi=\phi_1\oplus\phi_2:A[z;\sigma ]=Ag_1[zg_1;\sigma_1]\oplus
Ag_2[zg_2;\sigma_2]\longrightarrow {\cal S}\oplus M_2(\cue_5[x])$.
This isomorphism will allow us to make calculations in ${\cal
S}\oplus M_2(\cue_5[x])$ and then to reflect them in $A[z;\sigma ]$.
We are interested in the $S_3$-convolutional codes, these are
obtained by means of the direct summands left ideals of $A[z;\sigma
]$. Hence, we get the primitive idempotents of ${\cal S}$ and
$M_2(\cue_5[x])$, and then we apply $\phi^{-1}$ to them. Note that
it is easy to see that any idempotent in ${\cal S}$ or
$M_2(\cue_5[x])$ is primitive.

The idempotent matrices of ${\cal S}$ are of the form
$A=\left ( \begin{array}{cc} r & s \\ xt & 1-r \\
\end{array}\right )$ with $r(1-r)=xts$. First we suppose that
$r$ is different from  $0$ and $1$. We have two possibilities:
$x|r$ or $x|(1-r)$. If $x|r$, we call $r=kxd$, $s=dq$, $xt=kxp$,
$1-r=pq$. Then
$\left ( \begin{array}{cc} p & -d \\ xk & q \\
\end{array}\right )\cdot A=C$, where
$C=\left ( \begin{array}{cc} 0 & 0 \\ kx & q \\
\end{array}\right )$, and $\left ( \begin{array}{cc} p & -d \\ xk & q \\
\end{array}\right )$ has $\left ( \begin{array}{cc} q & d \\ -xk & p \\
\end{array}\right )$ as inverse in  ${\cal S}$.
Hence $^{\bullet}<A>=^{\bullet}<C>.$ If
$xk=\sum_{i=0}^n\alpha_ix^{i+1}$ and $q=\sum_{i=0}^m\beta_ix^i$,
then $\phi_1^{-1}(C)=\vae_2(\delta +\gamma )=u$ where
$\delta=\sum_{i=0}^n\alpha_iz^{2i+1}$, $\gamma
=\sum_{i=0}^m\beta_iz^{2i}$. Since $b_i\vae_2=\vae_2$ or
$4\vae_2$, we get a convolutional code of rank $1$, with the
$\cue_5[z]$-basis $\{ (\delta +\gamma,\delta +4\gamma ,\delta
+4\gamma ,\delta +4\gamma ,\delta +\gamma ,\delta +\gamma )\}$,
and complexity $max\{ 2deg(k)+1, 2deg(q)\}$.

In the second case, that is, when  $x| (1-r)$, we have
$^{\bullet}<A>=^{\bullet}<C>$, where  $C$ is now $C=\left ( \begin{array}{cc} k & q \\ 0 & 0 \\
\end{array}\right )$. In the same way as above, we get a convolutional code of rank
1, with basis $\{ (\delta +\gamma,\delta +4\gamma ,\delta +4\gamma
,\delta +4\gamma ,\delta +\gamma ,\delta +\gamma )\}$, and
complexity $max\{ 2deg(k), 2deg(q)+1\}$, where
$\delta=\sum_{i=0}^n\alpha_iz^{2i}$, $\gamma
=\sum_{i=0}^m\beta_iz^{2i+1}$.

Finally, we compute rank, basis, and complexity of the codes that
we get when $r=0,1$:

$\left ( \begin{array}{cc} 0 & s \\ 0 & 1 \\
\end{array}\right )$: rank 1, with basis $\{ (1,4,4,4,1,1) \}$, and complexity zero.
\medskip

$\left ( \begin{array}{cc} 1 & s \\ 0 & 0 \\
\end{array}\right )$: rank 1, with basis $\{ (1+\gamma ,4\gamma +1,4\gamma +1,4\gamma +1,1+\gamma ,1+ \gamma ) \}$,
and complexity $2deg(s)+1$. (If $s=\sum_{i=0}^n\alpha_ix^i$, then
$\gamma =\sum_{i=0}^n\alpha_iz^{2i+1}$).
\medskip

$\left ( \begin{array}{cc} 0 & 0 \\ xt & 1 \\
\end{array}\right )$: rank 1, with basis $\{ (1+\gamma ,4+\gamma ,4+\gamma ,4+\gamma ,1+\gamma ,1+\gamma ) \}$,
and complexity $2deg(t)+1$. (If $t=\sum_{i=0}^n\alpha_ix^i$, then
$\gamma =\sum_{i=0}^n\alpha_iz^{2i+1}$).
\medskip

$\left ( \begin{array}{cc} 1 & 0 \\ xt & 0 \\
\end{array}\right )$: rank 1, with basis $\{ (1,1 ,1 ,1,1,1 ) \}$,
and complexity zero. \bigskip

We resume all the above by stating that any minimal
$S_3$-convolutional code corresponding to an idempotent of ${\cal
S}$ has the  basis

\noindent $\{ (f(z),f(-z) ,f(-z) ,f(-z) ,f(z) ,f(z) )\}\;$ or
$\;\{ (f(z),-f(-z) ,-f(-z) ,$

\noindent$-f(-z) ,f(z) ,f(z) )\} ,$
 where $f(z)$ is a polynomial in $\cue_5[z]$,
$f(z)=\sum_{i=0}^na_iz^i$, with $\sum a_{2i}z^{2i}$ and $\sum
a_{2i+1}z^{2i+1}$ coprime (or, equivalently, $f(z)$ and $f(-z)$
coprime), or $\sum a_{2i}z^{2i}=1$ and $\sum a_{2i+1}z^{2i+1}=0$.
Hence the complexity is always $deg(f)$. In both cases, these
codes can be seen as codes of length $2$ by concatenation.

For several small $deg(f)$ we can compute the free distance of some
of these codes. For example, if $f(z)=bz+a$ with $a,b\neq 0$ the
code generated by $(f(z),f(-z) ,f(-z) ,f(-z) ,f(z) ,f(z) )$ has free
distance $12$ and so is a MDS code. It is also easy to see that if
$f(z)=a+bz+cz^2$ with $a,b,c\neq 0$, then the code generated by
$(f(z),f(-z) ,f(-z) ,f(-z) ,f(z) ,f(z) )$ has free distance $18$ and
so is a MDS code too.
\bigskip\bigskip

Now we focus our attention into the idempotents of
$M_2(\cue_5[x])$. Set $d=\phi_2^{-1}(B)=4(1\, 2\, 3)+(1\, 3\,
2)\in Ag_2$. Note that $d^2=2g_2$, hence $d^{2t}=2^tg_2$ and
$d^{2t+1}=2^tdg_2$.

We consider an idempotent matrix in
$M_2(\cue_5[x])$: $\left ( \begin{array}{cc} r & s \\ t & 1-r \\
\end{array}\right )$ with $r\neq 0,1$. Since $r(1-r)=ts$, we call $r=ah$, $s=hc$,
$t=ab$, $1-r=bc$. Then we have the following equalities:
$$\left ( \begin{array}{cc} a & c \\ -b & h \\
\end{array}\right )\cdot \left ( \begin{array}{cc} r & s \\ t & 1-r \\
\end{array}\right )\cdot \left ( \begin{array}{cc} h & -c \\ b & a \\
\end{array}\right ) =\left ( \begin{array}{cc} 1 & 0 \\ 0 & 0 \\
\end{array}\right ) ,$$
where $$\left ( \begin{array}{cc} a & c \\ -b & h \\
\end{array}\right )^{-1}=\left ( \begin{array}{cc} h & -c \\ b & a \\
\end{array}\right ) .$$
Then
$$\left ( \begin{array}{cc} a & c \\ -b & h \\
\end{array}\right )\cdot \left ( \begin{array}{cc} r & s \\ t & 1-r \\
\end{array}\right )=\left ( \begin{array}{cc} a & c \\ 0 & 0 \\
\end{array}\right ) .$$
Hence the left ideals generated by  $\left ( \begin{array}{cc} r & s \\ t & 1-r \\
\end{array}\right )$ and $\left ( \begin{array}{cc} a & c \\ 0 & 0 \\
\end{array}\right )$ are the same. So we only have to transform  $\left ( \begin{array}{cc} a & c \\ 0 & 0 \\
\end{array}\right )$ into an element of $A[z;\sigma ]$ and then calculate the associated convolutional code.

Let $a=\sum_{i=0}^n\alpha_ix^i, c=\sum_{i=0}^m\beta_ix^i\in
\cue_5[x]$. Then,

\begin{flushleft}
$
\left ( \begin{array}{cc} a & c \\ 0 & 0 \\
\end{array} \right ) =\left ( \begin{array}{cc} a & 0 \\ 0 & 0 \\
\end{array} \right ) +\left ( \begin{array}{cc} 0 & c \\ 0 & 0 \\
\end{array} \right ) =$
\end{flushleft}
\begin{flushright}
 $\sum_{i=0}^n\left ( \begin{array}{cc} x^i & 0 \\ 0 & x^i \\
\end{array} \right ) \left ( \begin{array}{cc} 1 & 0 \\ 0 & 0 \\
\end{array} \right ) \alpha_i  +\sum_{i=0}^m\left ( \begin{array}{cc} x^i & 0 \\ 0 & x^i \\
\end{array} \right ) \left ( \begin{array}{cc} 0 & 1 \\ 0 & 0 \\
\end{array} \right ) \beta_i.$
\end{flushright}

Hence,

$$\phi_2^{-1}\left ( \begin{array}{cc} a & c \\ 0 & 0 \\
\end{array} \right )
=\sum_{i=0}^n(zg_2)^id^i\vae_{3}\alpha_i+\sum_{i=0}^m(zg_2)^id^i\vae_{34}\beta_i$$
and so $$\phi^{-1}\left ( \begin{array}{cc} a & c \\ 0 & 0 \\
\end{array} \right )
=\sum_{i=0}^nz^id^i\vae_{3}\alpha_i+\sum_{i=0}^mz^id^i\vae_{34}\beta_i=u.$$(Note
that $g_2$ is the identity in $Ag_2$).

Set $a' =\sum_{i=0}^nz^i\alpha_id^i$, $c'
=\sum_{i=0}^mz^i\beta_id^i$. Breaking $a'$ and $c'$ according to
the parity of the $z$-degree of the monomials we write: $a_1 =\sum
z^{2i}\alpha_{2i}2^i$, $a_2=\sum z^{2i+1}2^i\alpha_{2i+1}$,$c_1
=\sum z^{2i}\beta_{2i}2^i$, $c_2=\sum z^{2i+1}2^i\beta_{2i+1}$.
Then $u=(a_1+da_2)\vae_3+(c_1+dc_2)\vae_{34}\in A[z;\sigma ]$.

In order to determinate the associated $S_3$-convolutional code, we
must calculate $b_i\vae_3$, $b_i\vae_{34}$,$b_i\vae_{3}d$ and
$b_i\vae_{34}d$, and then calculate $b_iu$. The final expression of
each $b_iu$ will be of the form
$b_iu=a_1u_{i1}+a_2u_{i2}+c_1u_{i3}+c_2u_{i4}$, with $u_{ij}\in A$.
This happens since $\sigma_2^2=I$. Taking this into account, with
the help of GAP software \cite{GAP4}, we get the generating matrix
whose files are the following:

$$\begin{array}{cl}
b_1u\mapsto w_1=&(2 a_1 + 3 a_2 + 4 c_2, 3 a_1 + 3 a_2 + 4 c_2, 4
a_2 +2 c_1 + 3 c_2, 2 a_1 + \\
 & 3 a_2 +3 c_1 + 3 c_2, 3 a_1 + 3 a_2 + 2 c_1 + 3 c_2,
  4 a_2 + 3 c_1 + 3 c_2),
\end{array}$$

$$\begin{array}{cl}
b_2u\mapsto w_2=& (3 a_1 + 2 a_2 + c_2, 2 a_1 + 2 a_2 + c_2, a_2 +
3 c_1 + 2 c_2,
  3 a_1 + 2 a_2 +  \\

& 2 c_1 +2 c_2, 2 a_1 + 2 a_2 + 3 c_1 + 2 c_2,
  a_2 + 2 c_1 + 2 c_2),
\end{array}$$

$$\begin{array}{cl}
b_3u\mapsto w_3=&  (a_2 + 2 c_1 + 2 c_2, 3 a_1 + 2 a_2 + 2 c_1 + 2
c_2, 2 a_1 + 2 a_2 + c_2, a_2 + \\
&3 c_1+2 c_2,3 a_1 + 2 a_2 + c_2, 2 a_1 + 2 a_2 + 3 c_1 + 2 c_2),
\end{array}$$

$$\begin{array}{cl}
b_4u\mapsto w_4=& (2 a_1 + 2 a_2 + 3 c_1 + 2 c_2, a_2 + 3 c_1 + 2
c_2, 3 a_1 + 2 a_2 + 2 c_1 +2 c_2,\\
&  2 a_1 + 2 a_2 + c_2, a_2 + 2 c_1 + 2 c_2,3 a_1 + 2 a_2 + c_2),
\end{array}$$

$$\begin{array}{cl}
 b_5u\mapsto w_5=& (4 a_2 + 3 c_1 + 3 c_2, 2 a_1 + 3 a_2 + 3 c_1 + 3 c_2,
  3 a_1 + 3 a_2 + 4 c_2, \\
&4 a_2 + 2 c_1 + 3 c_2,2 a_1 + 3 a_2 + 4 c_2,
  3 a_1 + 3 a_2 + 2 c_1 + 3 c_2),
\end{array}$$

$$\begin{array}{cl}
 b_6u\mapsto w_6=& (3 a_1 + 3 a_2 + 2 c_1 + 3 c_2, 4 a_2 + 2 c_1 +
3 c_2, 2 a_1 + 3 a_2 + 3 c_1 +\\
& 3 c_2,3 a_1 + 3 a_2 + 4 c_2, 4 a_2 + 3 c_1 + 3 c_2, 2 a_1 + 3
a_2 + 4 c_2) .
\end{array}$$

It is easy to see that $w_1=-w_2$, $w_4=-w_2-w_3$, $w_5=-w_3$ and
$w_6=w_2+w_3$. Therefore the code has rank  $2$, $\{ w_2,w_3\}$ is
a basis and the complexity is $max\{ 2deg(a), 2deg(c) \}$.

When the idempotent matrix of $M_2(\cue_5[x])$ has $r=0$ or $r=1$,
we can reduce its study to the above case. Concretely, we have

$$\left ( \begin{array}{cc} 0 & 1 \\ -1 & s \\
\end{array}\right )\cdot \left ( \begin{array}{cc} 0 & s \\ 0 & 1 \\
\end{array}\right )  =\left ( \begin{array}{cc} 0 & 1 \\ 0 & 0 \\
\end{array}\right ) ,$$

$$\left ( \begin{array}{cc} 0 & 1 \\ 1 & 0 \\
\end{array}\right )\cdot \left ( \begin{array}{cc} 0 & 0 \\ t & 1 \\
\end{array}\right )  =\left ( \begin{array}{cc} t & 1 \\ 0 & 0 \\
\end{array}\right ) ,$$

$$\left ( \begin{array}{cc} 1 & 0 \\ -t & 1 \\
\end{array}\right )\cdot \left ( \begin{array}{cc} 1 & 0 \\ t & 0 \\
\end{array}\right )  =\left ( \begin{array}{cc} 1 & 0 \\ 0 & 0 \\
\end{array}\right ) ,$$
where the left side matrices of the product are invertible in
$M_2(\cue_5[x])$, (the matrix $\left ( \begin{array}{cc} 1 & s \\ 0 & 0 \\
\end{array}\right )$ is not necessary to be reduced).

Therefore, all the minimal $S_3$-convolutional codes corresponding
to idempotents in the component $M_2(\cue_5[x])$ have the basis

\begin{flushleft}
$\{ (3 a_1 + 2 a_2 + c_2, 2 a_1 + 2 a_2 + c_2, a_2 + 3 c_1 + 2
c_2,
  3 a_1 + 2 a_2 + 2 c_1 + 2 c_2,$
\end{flushleft}
\begin{flushright}
  $2 a_1 + 2 a_2 + 3 c_1 + 2 c_2, a_2 + 2 c_1 + 2 c_2),$
\end{flushright}
\begin{flushleft}
  $(a_2 + 2 c_1 + 2 c_2, 3 a_1 + 2 a_2 + 2 c_1 + 2 c_2,
  2 a_1 + 2 a_2 + c_2, a_2 + 3 c_1 + 2 c_2,$
\end{flushleft}
\begin{flushright}
  $3 a_1 + 2 a_2 + c_2, 2 a_1 + 2 a_2 + 3 c_1 + 2 c_2) \}$
\end{flushright}
where $a_1 =\sum z^{2i}\alpha_{2i}2^i$, $a_2=\sum
z^{2i+1}2^i\alpha_{2i+1}$,$c_1 =\sum z^{2i}\beta_{2i}2^i$, $c_2=$

\noindent $\sum z^{2i+1}2^i\beta_{2i+1}$, and $a=\sum
\alpha_iz^i$, $c=\sum \beta_iz^i$ are any coprime polynomials in
$\cue_5[z]$, or $a=0$, $c=1$, or $a=1$, $c=0$. The rank is always
$2$ and the complexity is always $max\{ 2deg(a), 2deg(c)\}$. Note
that, in the above basis, the second vector is obtained from the
first one by permuting the components with $(1\, 5\, 6)(2\, 3\,
4)$.

\subsection{ The case of the identity permutation}

Now we study the $S_3$-convolutional codes that are obtained when
the automorphism maps $f_1$ to $f_1$. We can take, without lost of
generality, $\sigma =id_A$. Then

$A[z;\sigma ]=A[z]=Ag_1[z]\oplus Ag_2[z]=(A\vae_1[z]\oplus
A\vae_2[z])\oplus Ag_2[z]\cong (\cue_5[z]\oplus \cue_5[z])\oplus
M_2(\cue_5)[z]\cong (\cue_5[z]\oplus \cue_5[z])\oplus
M_2(\cue_5[z])$.

Hence $Ag_1[z]$ has only two idempotents different from $0$ and
$1$, concretely, $\vae_1$ and $\vae_2$, which generate two direct
summand left ideals of $Ag_1[z]$. The $S_3$-convolutional code
associated to $\vae_1$ has rank $1$, a basis is $\{ (1,1,1,1,1,1)
\}$, that is, it is a block code. The $S_3$-convolutional code
associated to $\vae_2$ has also rank $1$, a basis is $\{
(1,4,4,4,1,1) \}$, i.e., it is a block code too. These are the
only minimal codes to consider in the component $Ag_1[z]$.

Next, we study the component $Ag_2[z]$. In the same way that in
the case $\sigma \neq id_A$ above, we find idempotent elements in
$Ag_2[z]$ corresponding to the respective idempotent matrices in
$M_2(\cue_5[z])$.

We start with the same situation that in the case $\sigma \neq
id$. We consider an arbitrary idempotent matrix
$\left ( \begin{array}{cc} r & s \\ t & 1-r \\
\end{array}\right )$ with $r(1-r)=ts$ and $r\neq 0,1$.
We will reach to the same conclusion that in the case $\sigma \neq
id$: it is enough to work with the matrix
$\left ( \begin{array}{cc} a & c \\ 0 & 0 \\
\end{array}\right )$. Then, this matrix is performed into the element
$a\vae_3+c\vae_{34}$ of $Ag_2[z]$. The associated generating
matrices of the minimal $S_3$-convolutional codes  are obtained in a
similar way to the case $\sigma \neq id_A$: we only have to put in
those matrices $a_2=c_2=0$ and consider $a_1=a$, $c_1=c$ as
arbitrary coprime polynomials in $\cue_5[z]$. The generating matrix
of the code has the following rows:

$w_1=(2 a, 3 a, 2 c, 2 a + 3 c, 3 a + 2 c, 3 c)$,

$w_2=(3 a, 2 a, 3 c, 3 a + 2 c,
    2 a + 3 c, 2 c),$

    $w_3=(2 c, 3 a + 2 c, 2 a, 3 c, 3 a, 2 a + 3 c),$

    $w_4=(2 a + 3 c,
    3 c, 3 a + 2 c, 2 a, 2 c, 3 a)$

    $w_5=(3 c, 2 a + 3 c, 3 a, 2 c, 2 a,
    3 a + 2 c)$

    $w_6=(3 a + 2 c, 2 c, 2 a + 3 c, 3 a, 3 c, 2 a).$

Then  $$w_1=-w_2,\; w_4=-w_2-w_3,\; w_5=-w_3,\; w_6=w_2+w_3$$

Therefore the code has rank $2$, $\{ w_2,w_3\}$ is a basis and the
complexity is $max\{ 2deg(a), 2deg(c)\}$.

When $r=0$ or $r=1$, we can also reduce the matrices to reach out
the above case and then we get some particular cases.
\medskip

We can compute, by comparing column and row distances of the
generator matrices  and using GAP software \cite{GAP4}, all the
optimal minimal $(6,2,2)$ $S_3$-convolutional codes which are
obtained by means of the identity permutation. Note that the
Griesmer bound for the field $\cue_5$ and memory $m=1$ is
$G_5(6,2,2;1)=10$ which is less than the Singleton bound (which is
11). In the following table appears all the possible values for
$a$ and $c$ that produce non equivalent optimal codes in this
situation.

\begin{table}[h]
\label{tabla} {\small
$$\begin{array}{|c|cccccccc|}
\hline
a&z+1&z+1&z+1&z+1&z+1&z+1&z+2&z+2\\
c&2z+1&3z+1&4z+1&z+2&z+3&z+4&z+1&2z+2\\
\hline a&z+2&z+2&z+2&z+2 & z+3&z+3&z+3&z+3 \\ c& 3z+2&4z+2&z+3&z+4
&z+1&z+2&2z+3&3z+3 \\
\hline
a&z+3&z+3&z+4&z+4&z+4&z+4&z+4&z+4\\
c&4z+3&z+4&z+1&z+2&z+3&2z+4&3z+4&4z+4\\
\hline

\end{array}$$}
\caption{Optimal minimal $(6,2,2;1)_5$ $S_3$-convolutional codes
for the identity permutation}
\end{table}

\bigskip

\subsection{Weight-preserving ring automorphisms}

We will show that two different $\cue_5$-automorphisms  of $A$
that produce the same permutation on the set
 $\{ f_1,f_2,f_3\}$ also produce {\bf isometric} skew polynomial
 rings. This is a very important issue for guaranteing a complete
 classification of $S_3$-convolutional codes with controlled free
 distances into a concrete skew polynomial ring.

Any $\cue_5$-automorphism $\sigma$ of $A$ verifies $\sigma (g)=kh$
with  $k\in \cue_5-\{ 0\}$ and $g,h\in S_3$. More precisely,
$\sigma ((i\, j))=u \cdot (l\, m)$ with $u\in \{ 1, -1\}$ and
$\sigma ((i \, j\, m))=(j \, i\, m)$. Hence we get six
automorphisms for the case $\sigma (f_1)=f_2$ and six for the case
$\sigma (f_1)=f_1$.

Let $\sigma$, $\tau$ two $\cue_5$-automorphisms verifying $\sigma
(f_1)=f_2$ $\tau (f_1)=f_2$. We will define a ring isometry $\chi
:A[z; \sigma ]\rightarrow  A[z; \tau ]$. We have $$A[z; \sigma
]=(Ag_1)[zg_1; \sigma_1]\oplus (Ag_2)[zg_2; \sigma_2]$$ and $A[z;
\tau ]=(Ag_1)[zg_1; \tau_1]\oplus (Ag_2)[zg_2; \tau_2]$, where
$\sigma_i$, $\tau_i$ are the restriction automorphisms to $Ag_i$,
$i=1,2$.

The ring $Ag_1$ is generated as, $\cue_5$-vector space, by
$c_1=I+(1\, 2\, 3)+(1\, 3\, 2)$ and $c_2=(1\, 2)+(1\, 3)+(2\, 3)$.
Therefore $\sigma_1(c_1)=c_1=\tau_1(c_1)$ and
$\sigma_1(c_2)=-c_2=\tau_1(c_2)$. Let $\chi_1 : (Ag_1)[zg_1;
\sigma_1]\rightarrow (Ag_1)[zg_1; \tau_1]$ be simply the identity
map.

Now we will define a ring isometry $\chi_2 : (Ag_2)[zg_2;
\sigma_2]\rightarrow (Ag_2)[zg_2; \tau_2]$. Since $Ag_2\cong
M_2(\cue_5 )$ there is $u, v\in U(Ag_2)$ such that $\sigma_2
(ag_2)=u^{-1}ag_2u$ and $\tau_2 (ag_2)=v^{-1}ag_2v$, for all $a\in
A$. Let $\chi_2(zg_2)=zg_2v^{-1}u$ and $\chi_2|_{Ag_2}=id_{Ag_2}.$
Since $v^{-1}u$ is a unit in $Ag_2$, then $g_1+v^{-1}u$ is a unit
in $A$ and so $g_1+v^{-1}u=k\cdot g$ for some $k\in \cue_5-\{ 0\}$
and $g\in S_3$. Hence $zg_2v^{-1}u=zg_2(g_1+v^{-1}u)=(zg_2)kg$,
i.e., $\chi_2$ is weight-preserving. In order to see that $\chi_2$
is a ring isomorphism we only have to check that
$\chi_2(czg_2)=\chi_2(zg_2)\chi_2(\sigma_2(c))$ for all $c\in
Ag_2$. But,
$\chi_2(czg_2)=czg_2v^{-1}u=zg_2\tau_2(c)v^{-1}u=zg_2(v^{-1}cv)v^{-1}u=zg_2v^{-1}cu=zg_2v^{-1}u(u^{-1}cu)=\chi_2(zg_2)
\chi_2(\sigma_2(c))$.

Now it is clear that the sum $\chi =\chi_1+\chi_2:A[z;\sigma
]\rightarrow A[z; \tau ]$ is a well-defined ring isometry.

When $\sigma (f_1)=f_1$, $\tau (f_1)=f_1$ we have $$A[z;\sigma
]=(A\vae_1)[z\vae_1 ; \sigma_1]\oplus (A\vae_2)[z\vae_2 ;
\sigma_2]\oplus (A g_2)[zg_2 ; \sigma_3]$$ and $$A[z;\tau
]=(A\vae_1)[z\vae_1 ; \tau_1]\oplus (A\vae_2)[z\vae_2 ;
\tau_2]\oplus (A g_2)[zg_2 ; \tau_3]$$ where $\sigma_i$ and $\tau_i$
are the corresponding restriction automorphisms. However it is easy
to see that $\sigma_i=id_{A\vae_i}$, $\tau_i=id_{A\vae_i}$ for
$i=1,2$. Therefore, we simply take $\chi_i=id_{(A\vae_1)[z\vae_1 ]}$
for $i=1,2$. On the other hand, $Ag_2\cong M_2(\cue_5 )$ so we can
use the above idea to build an isometry $\chi_3 :(A g_2)[zg_2 ;
\sigma_3]\rightarrow (A g_2)[zg_2 ; \tau_3]$. Then $\chi
=\chi_1+\chi_2+\chi_3$ is the desired isometry.

\section{Conclusions} All the minimal $S_3$-convolutional
codes over $\cue_5$ have the parameters $(6,1,t )$ or $(6,2,2t )$
($t$ an arbitrary positive integer). If we compare this with the
parameters of minimal $\Z_6$-convolutional codes (that is,
$\sigma$-cyclic convolutional codes) we get the same result (see
\cite[Theorem 3.8]{OCCCL}). Hence all minimal group convolutional
codes of length $6$ over the field of five elements have
parameters $(6,1,t)$ or $(6,2,2t)$. The positive integer $t$
corresponds with the (constant) Forney indices of the code. Also
note that general group codes are significantly more complicated
than $\sigma$-cyclic convolutional ones. When $\sigma =id$, cyclic
convolutional codes are always block codes, however, this is not
the case for $S_3$-convolutional codes. Finally, some free
distances have been computed for these minimal $S_3$-convolutional
codes. The calculations show that MDS-convolutional codes (or
optimal codes) appear frequently in this setting. It would be
interesting to give some information on the free distance of group
convolutional codes in terms of the algebraic structure of the
groups.\bigskip

\centerline{\bf Acknowledgements } This work have been supported
by the grant BFM2002-02717 from DGES.

\end{document}